\documentclass[11pt,a4paper,reqno]{amsart}
\usepackage[margin=2.5cm]{geometry}
\usepackage{amssymb,mathtools}
\usepackage{enumitem}

\usepackage{xcolor}

\usepackage{hyperref}
\hypersetup{hypertexnames=false,colorlinks=true,linkcolor=blue,citecolor=red,pdftitle={Preventing cycles in Lotka--Volterra systems},pdfsubject={Lotka--Volterra systems: acyclicity and quantitative symmetric dominance}}

\theoremstyle{plain}
\newtheorem{theorem}{Theorem}[section]
\newtheorem{proposition}[theorem]{Proposition}
\newtheorem{lemma}[theorem]{Lemma}
\newtheorem{corollary}[theorem]{Corollary}

\theoremstyle{definition}
\newtheorem{definition}[theorem]{Definition}

\theoremstyle{remark}
\newtheorem{remark}[theorem]{Remark}

\newcommand{\R}{\mathbb R}
\newcommand{\one}{\mathbf1}
\DeclareMathOperator{\supp}{supp}
\newcommand{\norm}[1]{\lVert#1\rVert}
\newcommand{\diag}{\operatorname{diag}}
\newcommand{\Id}{\mathrm{Id}}
\newcommand{\cto}{\stackrel{c}{\longrightarrow}}
\newcommand{\wto}{\stackrel{w}{\longrightarrow}}

\numberwithin{equation}{section}

\title[Preventing cycles]{Preventing cycles in Lotka--Volterra systems
	}

\author{M.C. Bortolan, P. Kalita, J.A. Langa and R.O. Moura}

\address[M.C. Bortolan]{Departamento de Matemática, Centro de Ciências Físicas e Ma\-te\-má\-ti\-cas\\ Universidade Federal de Santa Catarina\\ Campus Florianópolis, Florianópolis, Brazil}
\email[M.C. Bortolan]{m.bortolan@ufsc.br}

\address[P. Kalita]{Faculty of Mathematics and Computer Science\\
Jagiellonian University\\
ul. \L{}ojasiewicza 6, 30-348 Krak\'ow, Poland}
\email[P. Kalita]{piotr.kalita@uj.edu.pl}

\address[J.A. Langa]{Departamento de Ecuaciones Diferenciales y Análisis Numérico\\
Universidad de Sevilla\\
Campus Reina Mercedes, 41012 Sevilla, Spain}
\email[J.A. Langa]{langa@us.es}

\address[R.O. Moura]{Artefact do Brasil Consultoria Estrat\'egica LTDA. \\
S\~ao Paulo, SP, Brazil}
\email[R.O. Moura]{rfomoura@gmail.com}

\thanks{(M.C. Bortolan) This study was financed in part by the Conselho Nacional de Desenvolvimento Cient\'{\i}fico e Tecnol\'ogico (CNPq, Brasil) under the Bolsa de Produtividade em Pesquisa 306136/2024-4 and by the Ministerio de Ciencia e Innovaci\'on, Agencia Estatal de Investigaci\'on (AEI) and FEDER (Espa\~na) under project PID2021-122991NB-C21.}
\thanks{(P. Kalita) This work was partially supported by the Spanish
Ministerio de Ciencia, Innovaci\'on y Universidades, Agencia
Estatal de Investigaci\'on (AEI), and FEDER under grant
PID2024-156228NB-I00, and by FAPESP, Brazil, under grant
2020/14075-6.}
\thanks{(J.A. Langa) This study was financed in part by the Ministerio de Ciencia, Innovaci\'on y Universidades, Agencia Estatal de Investigaci\'on, grant PID2024-156228NB-I00.}
\thanks{(R.O. Moura) This study was financed by Grants \#2022/04886-2 and \#2023/11798-5, S\~ao Paulo Research Foundation (FAPESP)}

\date{}



\begin{document}

\begin{abstract}
We study cycles of heteroclinic connections between equilibria of Lotka--Volterra systems with negative self-interactions. We prove that generalized strict diagonal dominance makes the invasion, connection, and weak graphs acyclic for every intrinsic growth vector. No hyperbolicity assumption is needed, and the associated semiflow is gradient with respect to its equilibria. For the broader class of Volterra--Lyapunov stable matrices, we show that two equilibria are joined by a heteroclinic connection exactly when the species that enter have positive invasion rates at the source and the species that leave have nonpositive invasion rates at the target. We show by examples that the three graphs can differ, that the dominance threshold is sharp, and that generalized dominance is not necessary for acyclicity. We also give a quantitative criterion for cycle exclusion in terms of the symmetric and skew parts of a weighted interaction matrix.
\end{abstract}

\maketitle

\section*{Declarations}

\textbf{Use of generative AI.} Generative AI tools (ChatGPT Pro, ChatGPT 6 Astra,
and Claude Opus~5) were used substantially in this research. The notion
of the weak graph (Definition~\ref{WG}), the characterization of heteroclinic
connections (Theorem~\ref{thm:exactcg}), and the quantitative criterion for
cycle exclusion (Theorem~\ref{thm:skew} and condition~\eqref{eq:uniform}) were
developed in dialogue with these tools. The key arguments in these parts
were proposed by AI in response to questions and conjectures formulated by
the authors, and were then refined and substantially simplified through
further iterations. For the proof of the result on  acyclicity via the scalar
potential $V_{\mathbf w}$ (Theorem~\ref{thm:potential}), the initial candidate Lyapunov function came from the authors, supported by numerical evidence from earlier work. Here AI tools supplied an initial proof and assisted with its subsequent refinement, verification, and presentation, as well as with drafting and editing the manuscript. The authors
formulated the research questions, selected and refined the arguments,
checked all statements, calculations, and proofs in detail, placed the
results in the context of the mathematical and ecological literature, and
rewrote the text for clarity and readability. The authors take full responsibility for the
manuscript and the correctness of its results. 

\medskip \textbf{Conflict of interest.} On behalf of all authors, the corresponding author states that there is no conflict of interest.

\medskip \textbf{Data availability.} No datasets were generated or analyzed; the work is theoretical.

\section{Introduction}\label{sec:1}

We consider the following Lotka--Volterra system for $N\ge1$ interacting species,
\begin{equation}\label{eq:lv}\tag{LV}
\dot{u}_i = u_i\Big(b_i + \sum_{j=1}^N a_{ij}u_j\Big) \quad \text{ for } i=1,\ldots,N,
\end{equation}
where $\mathbf{b}=(b_1,\ldots,b_N)^{\mathsf T}\in\mathbb{R}^N$ is the vector of \emph{intrinsic growth rates} and $A=(a_{ij})_{i,j=1}^N\in\mathbb{R}^{N\times N}$ is the \emph{interaction matrix}. The phase space is
\[
\mathbb{R}^N_{\ge 0} = \{\mathbf{u}=(u_1,\ldots,u_N)\in \mathbb{R}^N\colon u_i\ge 0 \text{ for } i=1,\ldots,N\}.
\]
The nonnegative orthant $\mathbb{R}^N_{\ge 0}$ and each of its coordinate faces are invariant for \eqref{eq:lv}. With each equilibrium we associate its support, that is the set of indexes, on which the variables are nonzero.  A heteroclinic trajectory that connects two equilibria with different supports is then a transition from one community to another.

We ask when heteroclinic connections between equilibria of \eqref{eq:lv} can form a directed cycle. Our main result excludes such cycles for all intrinsic growth rates when $A$ has negative diagonal entries and satisfies generalized strict diagonal dominance. This generalized strict diagonal dominance with negative diagonal implies Volterra--Lyapunov stability, but not conversely when $N\ge2$. A matrix $A$ is \emph{Volterra--Lyapunov stable} (\emph{VL stable}, for short) if there exists a positive diagonal matrix $H$ such that $HA+A^{\mathsf T}H\prec0$, that is, $HA+A^{\mathsf T}H$ is negative definite. We write
\[
\mathcal{V}
=\{A\in\mathbb{R}^{N\times N}\colon A\text{ is VL stable}\}.
\]
We recall the following standard consequence of Volterra--Lyapunov
stability, see \cite[Theorem~8]{AKLS}.
\begin{proposition}\label{prop.old}
Assume that $A\in\mathcal{V}$. For every $\mathbf{u}_0\in\mathbb{R}^N_{\ge0}$,
the problem \eqref{eq:lv} has a unique solution
$\mathbf{u}(\cdot,\mathbf{u}_0)\colon[0,\infty)\to\mathbb{R}^N_{\ge0}$
with $\mathbf{u}(0,\mathbf{u}_0)=\mathbf{u}_0$. The solution depends
continuously on time and initial data and defines a continuous semiflow
$\mathcal{T}=\{T(t):t\ge0\}$ by
$T(t)\mathbf{u}_0=\mathbf{u}(t,\mathbf{u}_0)$.
This semiflow possesses a global attractor $\mathcal{A}$.
\end{proposition}

Diagonal dominance means that the self-interaction of each species outweighs its interactions with the others. The matrix $A$ is strictly row diagonally dominant if in every row the diagonal entry exceeds in absolute value the sum of the absolute values of the other entries, that is, if
 \[
 |a_{ii}| > \sum_{j\ne i} |a_{ij}| \quad \text{for } i=1,\ldots,N.
 \]
In the generalized version the species carry positive weights, and we ask for a vector $\mathbf{v}>\mathbf{0}$ such that
\[
|a_{ii}|v_i>\sum_{j\ne i}|a_{ij}|v_j
\qquad\text{for every }i=1,\ldots,N.
\]
Equivalently, $A\operatorname{diag}(\mathbf{v})$ is strictly row diagonally dominant. One may equally weight the columns, since the same class of matrices is obtained by asking for a vector $\mathbf w>\mathbf 0$ such that $|a_{jj}|w_j>\sum_{i\ne j}|a_{ij}|w_i$ for every $j$, that is, $\operatorname{diag}(\mathbf w)A$ is strictly column diagonally dominant. 
Several classical conditions are equivalent to generalized strict diagonal dominance, we recall them in Proposition~\ref{th:equiv} of Appendix~\ref{sec:equiv}, see \cite{FiedlerPtak1962,Plemmons1977,BermanPlemmons1994,BCGM}. Here we use the one in terms of a spectral radius.
Throughout the paper we always assume $a_{ii}<0$ for $i=1,\ldots,N$. We define $B=(B_{ij})_{i,j=1}^N$ as the matrix given by
\begin{equation}\label{eq:B}
B_{ii}=0 \quad \text{ and } \quad B_{ij}=\frac{|a_{ij}|}{|a_{ii}|} \quad \text{ for } i\ne j.
\end{equation}
With
\[
\rho(B) = \max\{|\lambda|\colon \lambda \text{ is an eigenvalue of } B\},
\]
we define the class of matrices
\begin{equation}\label{def.H-}
	\mathcal{H}_-=\{A\in \mathbb{R}^{N\times N}\colon a_{ii}<0 \text{ for } i=1,\ldots,N \text{ and } \rho(B)<1\}.
\end{equation}

Diagonal dominance has appeared before in ecological stability theory. Hofbauer and So \cite[Theorem~2]{HofbauerSo2000} characterize global stability of an existing positive equilibrium for all admissible off-diagonal delays, and Virgo, Law and Emmerson \cite{VirgoLawEmmerson2006} use diagonal dominance in their study of sequentially assembled food webs. Stability of a positive equilibrium for positive initial data, however, says nothing about cycles between boundary equilibria, and these are what we study here.

In \cite[Corollary~4.3]{BKLM} we proved that strict row or column diagonal dominance excludes heteroclinic cycles in dimension three and cycles consisting only of single-species equilibria in arbitrary dimension, and we asked (Open Question~1 there) whether the restriction to single-species equilibria can be removed. Numerical evidence and a conjecture in this direction can be found in \cite{CallejaSolanas2026}. In Theorem~\ref{thm:main} we give a positive answer for strict dominance. The theorem is in fact more general, in that the dominance may be weighted, the supports of the equilibria and the growth vector are arbitrary, and hyperbolicity is not assumed. Under non-strict dominance heteroclinic cycles can occur, see Proposition~\ref{th:sharp}.

For any $\mathbf{x}\in\mathbb{R}^N_{\ge0}$, we write
$\supp\mathbf{x}=\{i:x_i>0\}$. At an equilibrium $\mathbf{x}$,
the \emph{invasion rate} of species $i$ is
\begin{equation}\label{eq:rates}
r_i(\mathbf{x})=b_i+(A\mathbf{x})_i = b_i + \sum_{j=1}^N a_{ij}x_j,
\end{equation}
which is zero for each resident species $i\in\supp\mathbf{x}$. For an absent species, it is its per-capita growth rate when introduced at low density.

For $A\in\mathcal V$, every nonempty principal submatrix is nonsingular.
Thus each support determines at most one equilibrium: its resident
coordinates solve $A_I\mathbf{x}_I=-\mathbf{b}_I$, where $A_I$ is the
principal submatrix indexed by $I$ and $\mathbf{x}_I,\mathbf{b}_I$ are
the corresponding subvectors. The empty support gives the zero equilibrium.

We use three directed graphs whose vertices are the nonnegative equilibria of \eqref{eq:lv}. The first is the equilibrium version of the \emph{invasion graph} ($\mathrm{IG}$) from \cite{Hofbauer}, in which we retain the strict sign inequalities. Replacing both strict inequalities by non-strict ones gives the \emph{weak graph} ($\mathrm{WG}$). The third is the \emph{connection graph} ($\mathrm{CG}$) of \cite{BKLM}.

\begin{definition}\label{WG}
Let $\mathbf{x}$ and $\mathbf{y}$ be distinct nonnegative equilibria of \eqref{eq:lv} with supports $I$ and $J$, respectively.
\begin{enumerate}[label={\bf (\alph*)}]
\item \textbf{Invasion graph (IG).} The edge $\mathbf{x} \longrightarrow \mathbf{y}$ exists if and only if $r_j(\mathbf{x}) > 0$ for every $j \in J \setminus I$, and $r_i(\mathbf{y}) < 0$ for every $i \in I \setminus J$.
\item \textbf{Weak graph (WG).} The edge $\mathbf{x}\wto \mathbf{y}$ exists if and only if
\begin{equation}\label{eq:weak}
r_j(\mathbf{x}) \ge 0  \text{ for every } j \in J \setminus I, \text{ and }  r_i(\mathbf{y}) \le 0 \text{ for every } i \in I \setminus J.
\end{equation}
\item \textbf{Connection graph (CG).} The edge $\mathbf{x}\cto\mathbf{y}$ exists if
and only if there is a bounded complete solution
$\boldsymbol{\xi}\colon\mathbb{R}\to\mathbb{R}^N_{\ge0}$ such that
\[
\mathbf{x}\stackrel{t \to -\infty}{\longleftarrow} \boldsymbol{\xi}(t) \stackrel{t \to \infty}{\longrightarrow} \mathbf{y}.
\]
\end{enumerate}
\end{definition}

In (a) and (b) we call the condition on the species in $J\setminus I$, which enter the community, the \emph{entry condition}, and the condition on the species in $I\setminus J$, which leave it, the \emph{exclusion condition}. Both are strict in $\mathrm{IG}$ and weak in $\mathrm{WG}$.

\medskip All three graphs have the same equilibrium vertices and contain no self-loops. When $A\in\mathcal V$, each support carries at most one equilibrium, so the vertices can also be indexed by the \emph{admissible} supports, or communities, that is, by the supports of equilibria. That there are no self-loops is only a convention. A nonconstant homoclinic solution would not show up in any of the graphs, and we exclude such solutions separately, under VL stability, in Proposition~\ref{lem:homoclinic}.

A weak cycle is a directed cycle in $\mathrm{WG}$. Its edges need not be realized by trajectories, but since every heteroclinic connection between distinct equilibria satisfies the weak inequalities (Proposition~\ref{lem:connection}), there are no heteroclinic cycles through distinct equilibria once $\mathrm{WG}$ is acyclic. Under VL stability we can say exactly which pairs of equilibria are connected. This is done in Section~\ref{sec:graphs}, where we also prove the inclusions $\mathrm{IG}\subseteq\mathrm{CG}\subseteq\mathrm{WG}$ and show by examples that either of them can be strict.

Generalized dominance is sufficient for the weak graph to be acyclic for every growth vector $\mathbf{b}$, but it is not necessary. For instance, it is known that the weak graph is also acyclic whenever $HA$ is symmetric and negative definite for some positive diagonal matrix $H$ (we recall this in Proposition~\ref{prop:symmetric}). In Section~\ref{sec:skew} we go further and show that it suffices that the skew-symmetric part of $HA$ be small, in appropriate sense, compared with its symmetric part.

Our main result, which we prove in Section~\ref{sec:dynamics}, is the following.

\begin{theorem}\label{thm:main}
Assume that $A\in\mathcal{H}_-$. For every $\mathbf{b}\in\mathbb{R}^N$, the following conclusions hold:
\begin{enumerate}[label={\bf (\alph*)}]
\item $A\in \mathcal{V}$;
\item all three graphs are acyclic, and \eqref{eq:lv} has no nonconstant homoclinic solution;
\item $\mathcal{T}$ has a global attractor $\mathcal{A}$;
\item $\mathcal{T}$ is gradient with respect to its equilibria: there is
a continuous function $\mathcal{L}\colon\mathbb{R}^N_{\ge0}\to\mathbb{R}$
such that $\mathcal{L}(T(t)\mathbf{u}_0)<\mathcal{L}(\mathbf{u}_0)$ for
every $t>0$ and every nonequilibrium initial state $\mathbf{u}_0$.
\end{enumerate}
Moreover, the $\alpha$-limit of each bounded complete solution of \eqref{eq:lv} is a single equilibrium.  
\end{theorem}

The proof rests on one inequality. Let $\mathbf x\wto\mathbf y$ be a weak edge. Since invasion rates vanish on supports, the exclusion condition in \eqref{eq:weak} allows us to bound the change $r_i(\mathbf y)-r_i(\mathbf x)=(A(\mathbf y-\mathbf x))_i$ of each invasion rate from above by an expression in which the off-diagonal entries of $A$ enter only through their absolute values (Corollary~\ref{cor:exclusion}). For column weights $\mathbf w$ of generalized dominance (see Proposition~\ref{th:equiv} of Appendix~\ref{sec:equiv}), the weighted sum of these bounds is negative, so that
\[
\sum_{i=1}^N w_ir_i(\mathbf y)<\sum_{i=1}^N w_ir_i(\mathbf x).
\]
The weighted sum of the invasion rates thus decreases strictly along every weak edge, and the weak graph has no cycle (Theorem~\ref{thm:potential} and Corollary~\ref{cor:acyclic}). In consequence connections are weak edges (Proposition~\ref{lem:connection}), VL stability excludes homoclinic solutions (Proposition~\ref{lem:homoclinic}), and the gradient structure then follows from \cite{BKLM}.

In community assembly, species arrive one at a time or in groups, and the basic question is whether the process ends in a
community that resists further invasions or continues indefinitely as species turnover~\cite{SchreiberRittenhouse2004,SpaakAdlerEllner2023,CallejaSolanas2026}. In particular, Spaak and Schreiber explicitly say in \cite{SpaakSchreiber2023} that ``identifying food-web topologies that ensure acyclic invasion graphs is an important open problem''. 
The order of arrival can decide which community is ultimately reached. This corresponds to the choice of the path in the graph of connections \cite{Fukami2015,SongFukamiSaavedra2021}.
Since there are finitely many possible equilibria, an assembly process that does not end must return to the one that was already visited, this means that there is a cycle in the connection graph~\cite{Song2025}.
The question if the graph is acyclic determines if the ecological assembly process will end, which also matters for ecological restoration
and the design of communities~\cite{ServanAllesina2021}.

The invasion graph is built from the resident communities and the signs of their invasion rates. If it is acyclic, and the other hypotheses of the invasion-graph framework hold, then permanence can be read off from these signs. The system is permanent if and only if every subcommunity, the empty one included, can be invaded by at least one absent species~\cite[Theorem~1]{Hofbauer} (see also \cite{SpaakSchreiber2023,ClarkEtAl2026,Song2025}). To check acyclicity one must find first all the system's equilibria and compute all invasion rates. The conditions on the model that guarantee acyclicity allows us to avoid this step to guarantee acyclicity~\cite{CallejaSolanas2026}.

We give such condition. It is possible to check in only based on the  interaction matrix without the need to construct the graph.  The diagonal values $a_{ii}$ must be negative, but condition allows both positive and negative off-diagonal values. Ordinary strict row or column diagonal dominance is the special case in which all weights are equal to one.

The paper is organized as follows. Section~\ref{sec:dynamics} contains the proof of Theorem~\ref{thm:main}, together with the construction of a potential that decreases along every edge of the weak graph. In Section~\ref{sec:graphs} we characterize the connections between equilibria without assuming hyperbolicity, and we show by examples that the threshold $\rho(B)=1$ is sharp and that the condition is not necessary for acyclicity. Quantitative criteria for cycle exclusion are given in Section~\ref{sec:skew}. The classical results on comparison matrices that we use are collected, with proofs, in Appendix~\ref{sec:equiv}.

\section{Acyclicity under generalized dominance}\label{sec:dynamics}

For a matrix $C=(c_{ij})_{i,j=1}^N\in\mathbb{R}^{N\times N}$ and nonempty
sets $I,J\subseteq\{1,\ldots,N\}$, write $C_{I,J}=(c_{ij})_{i\in I,j\in J}$
and $C_I=C_{I,I}$. For $\mathbf u\in\mathbb R^N$, write
$\mathbf u_I=(u_i)_{i\in I}$.

For $\mathbf{u}\in\R^N_{\ge0}$ we define $r(\mathbf{u})=(r_i(\mathbf{u}))_{i=1}^N=\mathbf{b}+A\mathbf{u}$, so that the $i$-th equation in the Lotka--Volterra system has the form $ \dot u_i=u_i r_i(\mathbf{u})$. At an equilibrium $\mathbf{x}\ge\mathbf{0}$ the numbers $r_i(\mathbf{x})$ are the invasion rates \eqref{eq:rates}, and $r_i(\mathbf{x})=0$ whenever $x_i>0$.
The \emph{comparison matrix} of $A$ is
\[
M=\diag(|a_{ii}|)\,(\Id-B),\qquad\text{that is,}\qquad M_{ii}=|a_{ii}|\quad\text{and}\quad M_{ij}=-|a_{ij}|\quad\text{for } i\ne j,
\]
where $B$ is given by \eqref{eq:B}. Here and below, $\Id$ denotes the identity matrix.
We begin with an elementary estimate, which compares $A\mathbf h$ with $-M|\mathbf h|$ coordinate by coordinate. It uses only that the diagonal of $A$ is negative.

\begin{lemma}\label{thm:comparison}
	Let $\mathbf{h}\in\R^N$ and fix $i\in\{1,\ldots,N\}$.
	If either $h_i\ge0$ or $(A\mathbf{h})_i\le0$, then
	\begin{equation}\label{eq:row}
		(A\mathbf{h})_i\le -(M|\mathbf{h}|)_i
		=-|a_{ii}|\,|h_i|+\sum_{j\ne i}|a_{ij}|\,|h_j|.
	\end{equation}
\end{lemma}

\begin{proof}
If $h_i\ge0$, then
	\[
	(A\mathbf{h})_i=-|a_{ii}|\,|h_i|+\sum_{j\ne i}a_{ij}h_j
	\le -|a_{ii}|\,|h_i|+\sum_{j\ne i}|a_{ij}|\,|h_j|.
	\]
	If $h_i<0$ and $(A\mathbf{h})_i\le0$, then $a_{ii}h_i=|a_{ii}|\,|h_i|$, since both factors are negative, and
	\[
	(A\mathbf{h})_i\le-(A\mathbf{h})_i
	=-|a_{ii}|\,|h_i|-\sum_{j\ne i}a_{ij}h_j
	\le -|a_{ii}|\,|h_i|+\sum_{j\ne i}|a_{ij}|\,|h_j|.
	\]
	This also proves \eqref{eq:row}.
\end{proof}

For the difference of two equilibria the lemma gives the following.
\begin{corollary}\label{cor:exclusion}
	Let $\mathbf{x},\mathbf{y}\ge\mathbf{0}$ be equilibria and assume that $r_i(\mathbf{y})\le0$ for every  $i\in\supp\mathbf{x}\setminus\supp\mathbf{y}$. Then
	\begin{equation}\label{eq:vector}
		r(\mathbf{y})-r(\mathbf{x})
		=A(\mathbf{y}-\mathbf{x})\le -M|\mathbf{y}-\mathbf{x}|.
	\end{equation}
\end{corollary}

\begin{proof}
	Put $\mathbf{h}=\mathbf{y}-\mathbf{x}$.
	For each coordinate, we verify one of the hypotheses of Lemma~\ref{thm:comparison}. There is nothing to prove when $h_i\ge0$. If $h_i<0$, then $x_i>0$ and $r_i(\mathbf{x})=0$.
	When $y_i>0$, we also have $r_i(\mathbf{y})=0$,
	when $y_i=0$, we must have $r_i(\mathbf{y})\le0$.
Thus, in the case $h_i<0$, we have $(A\mathbf{h})_i=r_i(\mathbf{y})-r_i(\mathbf{x})\le0$.
	Lemma~\ref{thm:comparison} therefore applies in every coordinate, and
	its estimate \eqref{eq:row} proves
	\eqref{eq:vector}.
\end{proof}

We next construct a potential $V_{\mathbf w}$ that orders the equilibria. It is built from a vector $\mathbf w$ of column weights of generalized diagonal dominance.
  
\begin{theorem}\label{thm:potential}
	Assume $A\in \mathcal{H}_-$ and let $\mathbf{w}>\mathbf{0}$ satisfy $M^{\mathsf T}\mathbf{w}>\mathbf{0}$, as in condition \textup{(iii)} of Proposition~\ref{th:equiv}. Put $\mathbf{c}=-A^{\mathsf T}\mathbf{w}$ and, on the set of nonnegative equilibria, define the potential
	\[
	V_{\mathbf{w}}(\mathbf{x})
	=\mathbf{w}^{\mathsf T}r(\mathbf{x}).
	\]
	Then $\mathbf{c}\ge M^{\mathsf T}\mathbf{w}>\mathbf{0}$ and $V_{\mathbf{w}}(\mathbf{x})=\mathbf{w}^{\mathsf T}\mathbf{b}-\mathbf{c}^{\mathsf T}\mathbf{x}$. Moreover,
	\begin{equation}\label{eq:potential}
		V_{\mathbf{w}}(\mathbf{y})-V_{\mathbf{w}}(\mathbf{x})
		=-\mathbf{c}^{\mathsf T}(\mathbf{y}-\mathbf{x})
		\le-(M^{\mathsf T}\mathbf{w})^{\mathsf T}|\mathbf{y}-\mathbf{x}|<0
	\end{equation}
	whenever $\mathbf{x}\ne\mathbf{y}$ are two equilibria such that $r_i(\mathbf{y})\le 0$ for every $i\in \supp \mathbf{x} \setminus \supp \mathbf{y}$. In particular, \eqref{eq:potential} holds along every edge $\mathbf{x}\wto\mathbf{y}$ of the weak graph. For $\mathbf{w}=M^{-\mathsf T}\mathbf{1}$ we have $\mathbf{c}\ge\mathbf{1}$, and the bound in \eqref{eq:potential} is $-\|\mathbf{y}-\mathbf{x}\|_1$.
\end{theorem}

\begin{proof}
	Since $\mathbf{w}\ge\mathbf{0}$, for every $j$ we have
	\[
	c_j=|a_{jj}|w_j-\sum_{i\ne j}a_{ij}w_i
	\ge |a_{jj}|w_j-\sum_{i\ne j}|a_{ij}|w_i=(M^{\mathsf T}\mathbf{w})_j>0.
	\]
	By $r(\mathbf{x})=\mathbf{b}+A\mathbf{x}$, we have $V_{\mathbf{w}}(\mathbf{x})=\mathbf{w}^{\mathsf T}\mathbf{b}+(A^{\mathsf T}\mathbf{w})^{\mathsf T}\mathbf{x}=\mathbf{w}^{\mathsf T}\mathbf{b}-\mathbf{c}^{\mathsf T}\mathbf{x}$, which also gives the equality in \eqref{eq:potential}.
	Let $\mathbf x\ne\mathbf y$ be as in the statement. Since $\mathbf w\ge\mathbf 0$, Corollary~\ref{cor:exclusion} gives
	\[
		V_{\mathbf{w}}(\mathbf{y})-V_{\mathbf{w}}(\mathbf{x})
	=\mathbf{w}^{\mathsf T}\bigl(r(\mathbf{y})-r(\mathbf{x})\bigr)  \le -\mathbf{w}^{\mathsf T}M|\mathbf{y}-\mathbf{x}|
		=-(M^{\mathsf T}\mathbf{w})^{\mathsf T}|\mathbf{y}-\mathbf{x}|.
	\]
	The right-hand side is negative, because $M^{\mathsf T}\mathbf{w}>\mathbf{0}$ and the nonnegative vector $|\mathbf{y}-\mathbf{x}|$ is not zero. Every weak edge satisfies the hypothesis on $\mathbf{x}$ and $\mathbf{y}$ by \eqref{eq:weak}. Finally, if $\mathbf{w}=M^{-\mathsf T}\mathbf{1}$, then $\mathbf{w}>\mathbf{0}$ and $M^{\mathsf T}\mathbf{w}=\mathbf{1}$ by Proposition~\ref{th:equiv}, so that $\mathbf{c}\ge\mathbf{1}$ and $(M^{\mathsf T}\mathbf{w})^{\mathsf T}|\mathbf{y}-\mathbf{x}|=\|\mathbf{y}-\mathbf{x}\|_1$.
\end{proof}

Using the above result we can exclude weak cycles. Recall that a weak cycle, that is, a directed cycle in the weak graph, consists of $m\ge2$ pairwise distinct nonnegative equilibria $\mathbf{x}_1,\ldots,\mathbf{x}_m$ with
\[
\mathbf{x}_1 \wto \mathbf{x}_2 \wto \cdots \wto \mathbf{x}_m \wto \mathbf{x}_1.
\]

\begin{corollary}\label{cor:acyclic}
If $A\in \mathcal{H}_-$, $\mathrm{WG}$ has no weak cycle.
\end{corollary}
\begin{proof}
Suppose such a cycle exists. Set $\mathbf{x}_{m+1}=\mathbf{x}_1$ and fix $\mathbf{w}>\mathbf{0}$ with $M^{\mathsf T}\mathbf{w}>\mathbf{0}$.  Theorem~\ref{thm:potential} gives a strict decrease
of $V_{\mathbf{w}}$ at each step, which is impossible around a cycle.
\end{proof}

\begin{remark}\label{rem:local}
Note that even if $A\not\in \mathcal{H}_-$, if, for some subcommunity $I\subset \{1,\ldots,n\}$, the minor $A_I$ is in $\mathcal{H}_-$, then the above argument excluded cycles between the communities being subsets of $I$.    
\end{remark}

\begin{remark}
	The weighted invasion rates that we use in our cycle exclusion criterion are linked with products $P_{w}(u)=\Pi_{i=1}^N u_i^{w_i}$, which are used is classical  criteria for permanence (see \cite{Hofbauer1981,  Jansen1987}). Along strictly positive solutions we have
	\[
	\frac{d}{dt}\log P_{\mathbf w}(\mathbf u(t))
	=\mathbf w^{\mathsf T}r(\mathbf u(t))=V_{\mathbf w}(\mathbf u(t)).
	\]
	Although $P_{\mathbf w}$ vanishes on the boundary and its logarithm does not make sense there, this logarithmic growth rate extends continuously there and equals to our function $V_{\mathbf w}$ at equilibria. Classical permanence criteria use positivity of suitable averages of this growth rate while in our result this same function decreases along edges. 
\end{remark}

Each coordinate of a complete solution $\boldsymbol\xi$ of \eqref{eq:lv} is either identically zero or positive for all times, so the support $\supp\boldsymbol\xi(t)$ does not depend on $t\in\mathbb R$. We denote this support by $\supp\boldsymbol\xi$. The following elementary estimate, for a coordinate that vanishes in the limit, will be used several times.

\begin{lemma}\label{lem:rates}
Let $\boldsymbol{\xi}$ be a complete solution of \eqref{eq:lv}, let $\mathbf{x}$ be an equilibrium, and let $j$ be an index such that $\xi_j$ is positive and $x_j=0$. Let $\varepsilon>0$.
\begin{enumerate}[label={\bf (\alph*)}]
\item If $\boldsymbol{\xi}(t)\to\mathbf{x}$ as $t\to-\infty$, then there is $T_\varepsilon\in\mathbb{R}$ such that, for every $t_0\le T_\varepsilon$,
\[
\xi_j(t_0)e^{(r_j(\mathbf{x})+\varepsilon)(t-t_0)}\le\xi_j(t)\le\xi_j(t_0)e^{(r_j(\mathbf{x})-\varepsilon)(t-t_0)}
\qquad\text{for all } t\le t_0.
\]
Consequently, $r_j(\mathbf{x})\ge0$.
\item If $\boldsymbol{\xi}(t)\to\mathbf{x}$ as $t\to\infty$, then there is $T_\varepsilon\in\mathbb{R}$ such that, for every $t_0\ge T_\varepsilon$,
\[
\xi_j(t_0)e^{(r_j(\mathbf{x})-\varepsilon)(t-t_0)}\le\xi_j(t)\le\xi_j(t_0)e^{(r_j(\mathbf{x})+\varepsilon)(t-t_0)}
\qquad\text{for all } t\ge t_0.
\]
Consequently, $r_j(\mathbf{x})\le0$.
\end{enumerate}
\end{lemma}
\begin{proof}
We have
\[
\frac{d}{dt} \log \xi_j(t) = \frac{\xi'_j(t)}{\xi_j(t)} = b_j + \sum_{k=1}^N a_{jk} \xi_k(t) \to r_j(\mathbf{x}),
\]
as $t\to-\infty$ in case (a) and as $t\to\infty$ in case (b).
In case (a), choose $T_\varepsilon$ such that $\bigl|\frac{d}{dt}\log\xi_j(t)-r_j(\mathbf{x})\bigr|\le\varepsilon$ for all $t\le T_\varepsilon$. For $t\le t_0\le T_\varepsilon$, integration over $[t,t_0]$ gives
\[
(r_j(\mathbf{x})-\varepsilon)(t_0-t)\le\log\xi_j(t_0)-\log\xi_j(t)\le(r_j(\mathbf{x})+\varepsilon)(t_0-t),
\]
which is the asserted estimate. If $r_j(\mathbf{x})<0$, take $\varepsilon=-r_j(\mathbf{x})/2$. The lower bound then tends to infinity as $t\to-\infty$, contradicting $\xi_j(t)\to x_j=0$. Thus $r_j(\mathbf{x})\ge0$.

Case (b) is proved in the same way, by integration over $[t_0,t]$. If $r_j(\mathbf{x})>0$, the lower bound with $\varepsilon=r_j(\mathbf{x})/2$ tends to infinity as $t\to\infty$, again contradicting $\xi_j(t)\to0$.
\end{proof}

It follows that every connection between distinct equilibria is an edge of the weak graph. Neither VL stability nor hyperbolicity of equilibria is needed for this.
\begin{proposition}\label{lem:connection}
Let $\mathbf{x}$ and $\mathbf{y}$ be two distinct equilibria of \eqref{eq:lv}. If $\mathbf{x}\cto \mathbf{y}$ then $\mathbf{x}\wto \mathbf{y}$.
\end{proposition}
\begin{proof}
Set $I=\supp\mathbf{x}$ and $J=\supp\mathbf{y}$, and let $\boldsymbol{\xi}$ be the bounded complete solution of \eqref{eq:lv} for which
\[
\mathbf{x} \stackrel{t\to -\infty}{\longleftarrow} \boldsymbol{\xi}(t) \stackrel{t\to \infty}{\longrightarrow} \mathbf{y}.
\]
For $j\in J\setminus I$, the coordinate $\xi_j$ is positive, because $y_j>0$, and $x_j=0$. Lemma~\ref{lem:rates}(a) therefore gives $r_j(\mathbf{x})\ge0$. For $i\in I\setminus J$, the coordinate $\xi_i$ is positive, because $x_i>0$, and $y_i=0$. Lemma~\ref{lem:rates}(b), applied to the forward limit $\mathbf{y}$, gives $r_i(\mathbf{y})\le0$.
\end{proof}

We compare a solution with an equilibrium on the boundary by means of a Lyapunov function of Goh type \cite{Goh1977, Takeuchi1996}. In the next lemma we compute its derivative, without using VL stability.

\begin{lemma}\label{lem:goh}
Let $H=\diag(h_1,\ldots,h_N)$ be a positive diagonal matrix, let $\mathbf{x}$ be an equilibrium of \eqref{eq:lv} with support $I$, and let $\boldsymbol{\xi}$ be a complete solution of \eqref{eq:lv} whose support $S$ contains $I$. For a set $P$ with $I\subseteq P\subseteq S$, define
\[
L_{\mathbf{x},P}(t)=\sum_{i\in I}h_i\left(\xi_i(t)-x_i-x_i\log\frac{\xi_i(t)}{x_i}\right)
       +\sum_{j\in P\setminus I}h_j \xi_j(t)
\]
and $\mathbf{q}(t)=\boldsymbol{\xi}_P(t)-\mathbf{x}_P$. Then $L_{\mathbf{x},P}(t)\ge0$ and
\begin{equation}\label{eq:goh}
\frac{d}{dt}L_{\mathbf{x},P}(t)
=\sum_{j\in P\setminus I}h_j\xi_j(t)r_j(\mathbf{x})
+\tfrac12\mathbf{q}(t)^{\mathsf T}(H_PA_P+A_P^{\mathsf T}H_P)\mathbf{q}(t)
+\mathbf{q}(t)^{\mathsf T}H_PA_{P,S\setminus P}\,\boldsymbol{\xi}_{S\setminus P}(t),
\end{equation}
where the last term is absent if $P=S$. If $\boldsymbol{\xi}(t)\to\mathbf{x}$ as $t\to-\infty$ or as $t\to\infty$, then $L_{\mathbf{x},P}(t)\to0$ in the same limit.
\end{lemma}
\begin{proof}
Each summand of $L_{\mathbf{x},P}$ is nonnegative, because $s-x_i-x_i\log(s/x_i)\ge0$ for $s>0$, and each summand tends to zero when $\boldsymbol{\xi}(t)\to\mathbf{x}$. Since $\dot\xi_i=\xi_ir_i(\boldsymbol{\xi})$, and $x_j=0$ for $j\in P\setminus I$,
\[
\frac{d}{dt}L_{\mathbf{x},P}(t)=\sum_{i\in P}h_i\bigl(\xi_i(t)-x_i\bigr)r_i(\boldsymbol{\xi}(t)).
\]
Moreover, $r_i(\boldsymbol{\xi}(t))=r_i(\mathbf{x})+\sum_{k\in S}a_{ik}(\xi_k(t)-x_k)$, where $r_i(\mathbf{x})=0$ for $i\in I$ and $x_k=0$ for $k\in S\setminus P$. Splitting the sum over $k$ into the parts with $k\in P$ and $k\in S\setminus P$ gives \eqref{eq:goh}, because $\mathbf{q}^{\mathsf T}H_PA_P\mathbf{q}=\tfrac12\mathbf{q}^{\mathsf T}(H_PA_P+A_P^{\mathsf T}H_P)\mathbf{q}$.
\end{proof}

A \emph{nonconstant homoclinic solution} is a nonconstant bounded complete solution $\boldsymbol{\xi}$ with $\boldsymbol{\xi}(t)\to\mathbf{x}$ as $t\to-\infty$ and as $t\to\infty$ for the same equilibrium $\mathbf{x}$. The connection graph has no loops and does not see such a solution, but the gradient criterion of \cite[Corollary~3.15]{BKLM} counts it as a cycle, so we have to exclude it separately. For this, VL stability is enough.

\begin{proposition}\label{lem:homoclinic}
If $A\in \mathcal{V}$, \eqref{eq:lv} has no nonconstant homoclinic solution.
\end{proposition}
\begin{proof}
Choose $H=\diag(h_1,\ldots,h_N)>0$ with $HA+A^{\mathsf T}H\prec0$.
Suppose $\boldsymbol{\xi}$ is a nonconstant homoclinic solution with limit
$\mathbf{x}$. Let $S=\supp\boldsymbol{\xi}$ and put $I=\supp\mathbf{x}\subseteq S$.
For $j\in S\setminus I$, the coordinate $\xi_j$ is positive and tends to
zero as $t\to\infty$. Lemma~\ref{lem:rates}(b)
therefore gives $r_j(\mathbf{x})\le0$.

Let $L_{\mathbf{x},S}$ be the function of Lemma~\ref{lem:goh} for $P=S$. In \eqref{eq:goh} the last term is then absent, the first sum is nonpositive, and the quadratic term is strictly negative, because $H_SA_S+A_S^{\mathsf T}H_S\prec0$ and uniqueness of solutions implies $\boldsymbol{\xi}(t)\ne\mathbf{x}$ for every $t$. Thus $L_{\mathbf{x},S}$ is strictly decreasing. Since it is nonnegative and tends to zero as $t\to-\infty$, this is impossible.
\end{proof}

For VL stable matrices, the dynamics is gradient as soon as the weak graph is acyclic.
\begin{corollary}\label{cor:wggradient}
If $A\in \mathcal{V}$ and $\mathrm{WG}$ is acyclic, then the semiflow $\mathcal{T}$ associated with \eqref{eq:lv} in $\mathbb{R}^N_{\ge 0}$ is gradient with respect to its equilibria, and the $\alpha$-limit of each bounded complete solution of \eqref{eq:lv} is a single equilibrium.
\end{corollary}
\begin{proof}
By Proposition~\ref{prop.old}, $\mathcal{T}$ has a global attractor. As $A\in \mathcal{V}$, every principal submatrix $A_I$ is nonsingular. Indeed, if $H$ is a positive diagonal matrix for which $HA+A^{\mathsf T}H\prec 0$, then $H_IA_I+A_I^{\mathsf T}H_I\prec 0$, and if $A_I\mathbf{u}=\mathbf 0$ for some $\mathbf{u}\ne \mathbf{0}$ we obtain $\mathbf{u}^{\mathsf T}(H_IA_I+A_I^{\mathsf T}H_I)\mathbf{u}=0$. Therefore, there are at most $2^N$ equilibria, one per support.

By Proposition~\ref{lem:connection}, every heteroclinic cycle through distinct equilibria is a weak cycle, so $\mathrm{CG}$ is acyclic, and by Proposition~\ref{lem:homoclinic} there are no nonconstant homoclinic solutions. By VL stability, $\omega$-limits are always single equilibria. The gradient and $\alpha$-limit conclusions follow
from \cite[Corollaries 3.15 and 3.16]{BKLM}, using the former for $N\ge3$ and the latter
for $N\le2$.
\end{proof}

\begin{proof}[Proof of Theorem~\ref{thm:main}]
We prove the classical implication $A\in\mathcal H_-\Rightarrow A\in\mathcal V$,
with an explicit $H$, in
Proposition~\ref{lem:vl} of Appendix~\ref{sec:equiv}. This is (a), and (c) then follows from Proposition~\ref{prop.old}.

By Corollary~\ref{cor:acyclic}, $\mathrm{WG}$ is acyclic. So is $\mathrm{IG}$, which is a subgraph of $\mathrm{WG}$ by definition, and so is $\mathrm{CG}$, by Proposition~\ref{lem:connection}. Since Proposition~\ref{lem:homoclinic} excludes homoclinic solutions, this proves (b). Finally, (d) and the assertion on $\alpha$-limits follow from Corollary~\ref{cor:wggradient}.
\end{proof}

The scalar potential from  Theorem~\ref{thm:potential} also has an ecological interpretation.
At an equilibrium $\mathbf{x}$ with support $I$, the invasion
rates of the resident species vanish, so
\[
V_{\mathbf{w}}(\mathbf{x})
=\mathbf{w}^{\mathsf T}r(\mathbf{x})
=\sum_{i\notin I}w_i r_i(\mathbf{x}).
\]
Thus $V_{\mathbf{w}}(\mathbf{x})$ is a sum, with positive weights, of the invasion rates of the
absent species. By
Theorem~\ref{thm:potential}, it decreases strictly along every edge of
the weak graph, and hence along every edge of the invasion and connection
graphs. Only the weighted sum decreases; the invasion rate of a single absent species may increase.

By Theorem~\ref{thm:potential}, the same ordering is expressed by the weighted equilibrium abundance $\mathbf{c}^{\mathsf T}\mathbf{x}=\mathbf{w}^{\mathsf T}\mathbf{b}-V_{\mathbf{w}}(\mathbf{x})$. Its weights $\mathbf{c}$ are positive, and it increases by at least $(M^{\mathsf T}\mathbf{w})^{\mathsf T}|\mathbf{y}-\mathbf{x}|$ along every edge $\mathbf{x}\wto\mathbf{y}$ of the weak graph.
If the variables represent biomass densities,
this quantity can be interpreted as a weighted equilibrium biomass.

\section{The three graphs and non-hyperbolic connections}\label{sec:graphs}

By Proposition~\ref{lem:connection}, every connection between distinct
equilibria is a weak edge. We now determine when a weak edge is realized by a
connection.

For $A\in\mathcal V$, each admissible support determines
a unique equilibrium, so the three graphs have a common finite vertex set;
when all equilibria are hyperbolic, the graphs coincide by \cite{BKLM}.
We prove a connection criterion that also covers non-hyperbolic equilibria.

\subsection{Connections without hyperbolicity}
We first prove that every species present along a complete solution but absent at its backward limit has a strictly positive invasion rate there.

\begin{proposition}\label{lem:strictentry}
Let $A\in \mathcal{V}$ and let $\boldsymbol{\xi}$ be a complete solution of \eqref{eq:lv} with $\boldsymbol{\xi}(t)\to \mathbf{x}$ as $t\to-\infty$, for some equilibrium $\mathbf{x}$ of \eqref{eq:lv}. Then $r_j(\mathbf{x})>0$ for every $j\in \supp \boldsymbol \xi \setminus \supp \mathbf{x}$.
\end{proposition}
\begin{proof}
Choose a positive diagonal matrix $H=\diag(h_1,\ldots,h_N)$ with $HA+A^{\mathsf T}H\prec0$. Let $S=\supp \, \boldsymbol\xi$ and $I=\supp \mathbf{x}\subseteq S$. By Lemma~\ref{lem:rates}(a), $r_j(\mathbf{x})\ge0$ for every $j\in S\setminus I$. We consider
\[
Z_1=\{j\in S\setminus I\colon r_j(\mathbf{x})>0\}\quad \text{and}\quad  Z_2=\{j\in S\setminus I\colon r_j(\mathbf{x})=0\}.
\]
The sets $I,Z_1,Z_2$ partition $S$. Suppose $Z_2\ne\varnothing$. Let $L_{\mathbf{x},P}$ be the function of Lemma~\ref{lem:goh} for $P=I\cup Z_2$, so that $S\setminus P=Z_1$. It is nonnegative, and $\ell(t)\to0$ as $t\to-\infty$. Since $r_j(\mathbf{x})=0$ for $j\in Z_2$, the first sum in \eqref{eq:goh} vanishes. With $\mathbf{q}(t)=\boldsymbol{\xi}_P(t)-\mathbf{x}_P$ we therefore obtain
\[
\tfrac{d}{dt}L_{\mathbf{x},P}(t)=\tfrac12\mathbf q^{\mathsf T}(t)(H_PA_P+A_P^{\mathsf T}H_P)\mathbf q(t)+\mathbf q^{\mathsf T}(t)H_PA_{P,Z_1}\boldsymbol\xi_{Z_1}(t)\le-C_0\|\mathbf q(t)\|_2^2+C_1\|\boldsymbol\xi_{Z_1}(t)\|_2^2,
\]
for constants $C_0>0$ and $C_1\ge0$, by VL stability and Young's inequality. Here $\|\cdot\|_2$ is the Euclidean norm.

Put $\alpha=\tfrac12\min_{j\in Z_1}r_j(\mathbf{x})>0$ if $Z_1\ne\varnothing$, and $\alpha=1$ if $Z_1=\varnothing$. Using Lemma~\ref{lem:rates}(a) with $\varepsilon=\alpha$ we obtain that for some $C_2>0$ and $t_0\in \mathbb{R}$ we have   $\xi_j(t)\le C_2e^{\alpha t}$ for every $j\in Z_1$ and $t\le t_0$, because $r_j(\mathbf{x})-\alpha\ge\alpha$ for $j\in Z_1$. If $Z_1=\varnothing$, we can take any $C_2$ and $t_0$. In both cases the estimate for $\tfrac{d}{dt}L_{\mathbf{x},P}$ gives
\[
\tfrac{d}{dt}L_{\mathbf{x},P}(t) \le C_3e^{2\alpha t}  \quad \text{ for every } t\le t_0,
\]
where $C_3=C_1C_2^2\,\#Z_1\ge0$.  Integration from $-\infty$ to $t$, $t<t_0$, yields
\[
 0\le L_{\mathbf{x},P}(t)\le \tfrac{C_3}{2\alpha} e^{2\alpha t}.
\]
For $j\in Z_2$, this implies $\xi_j(t)\le \tfrac{C_3}{2\alpha h_j} e^{2\alpha t}$ for all $t<t_0$. But $r_j(\mathbf{x})=0$ for $j\in Z_2$, so Lemma~\ref{lem:rates}(a) with $\varepsilon=\alpha$ gives $T<t_0$ such that 
\[
\xi_j(t)  \ge \xi_j(T)e^{\alpha(t-T)}\quad \text{for every}\quad t\le T.
\]
Therefore, for every $t\le T$ we have
\[
0 < \xi_j(T)e^{-\alpha T} \le \xi_j(t)e^{-\alpha t} \le   \tfrac{C_3}{2\alpha h_j} e^{\alpha t}, \]
which gives a contradiction when $t\to-\infty$. Thus $Z_2$ is empty.
\end{proof}
By the preceding proposition, the strict entry condition is necessary for a connection. We now prove that, together with the weak exclusion condition, it is also sufficient. In the proof we use two results of \cite{AKLS}: Theorem~5, on the global stability of the equilibrium that no absent species can invade, and Theorem~18, which provides the connecting solution.
\begin{theorem}\label{thm:exactcg}
Let $A\in \mathcal{V}$. For two distinct nonnegative equilibria $\mathbf{x},\mathbf{y}$ of \eqref{eq:lv} with supports $I,J$, respectively, we have $\mathbf{x}\cto \mathbf{y}$ if and only if
\begin{equation}\label{eq:mixed}
 r_j(\mathbf{x})>0 \text{ for } j\in J\setminus I, \quad \text{ and } \quad r_i(\mathbf{y})\le0 \text{ for } i\in I\setminus J.
\end{equation}
Consequently,
\begin{equation}\label{eq:graphs}
\mathrm{IG}\subseteq\mathrm{CG}\subseteq\mathrm{WG}.
\end{equation}
If the invasion rates of all absent species at all equilibria are nonzero (that is, all equilibria are hyperbolic), then all three
graphs coincide.
\end{theorem}
\begin{proof}
Suppose first that $\mathbf{x}\cto \mathbf{y}$. Then every species in $J\setminus I$ is positive along a complete solution $\boldsymbol{\xi}$ connecting $\mathbf{x}$ and $\mathbf{y}$. Proposition~\ref{lem:strictentry} gives the strict entry condition, that is, $r_j(\mathbf{x})>0$ for all $j\in J\setminus I$.  Proposition~\ref{lem:connection} gives the weak exclusion condition, that is, $r_i(\mathbf{y})\le 0$ for $i\in I\setminus J$.

Conversely, assume \eqref{eq:mixed} and put $U=I\cup J$. Consider the subsystem of \eqref{eq:lv} formed by the species in $U$, with interaction matrix $A_U$ and growth vector $\mathbf b_U$. Its solutions are the solutions of \eqref{eq:lv} whose coordinates outside $U$ vanish, and $A_U$ is VL stable: if $H$ is a positive diagonal matrix with $HA+A^{\mathsf T}H\prec0$, then $H_UA_U+A_U^{\mathsf T}H_U\prec0$.

The equilibrium $\mathbf y$ is \emph{saturated} for this subsystem: $r_j(\mathbf y)=0$ for $j\in J$, and $r_i(\mathbf y)\le0$ for $i\in U\setminus J=I\setminus J$ by \eqref{eq:mixed}. Equivalently, $\mathbf y_U$ solves the linear complementarity problem
\[
\begin{cases}&\mathbf u\ge\mathbf 0\qquad \textrm{and}\qquad \mathbf b_U+A_U\mathbf u\le\mathbf 0,\\
	& \mathbf u^{\mathsf T}(\mathbf b_U+A_U\mathbf u)=0 .
	\end{cases}
\]
By \cite[Theorem~5]{AKLS}, which restates \cite[Theorem~3.2.1]{Takeuchi1996}, this problem has exactly one solution when $A_U$ is VL stable, and the corresponding equilibrium attracts every solution of the subsystem whose coordinates in its support are positive. Hence $\mathbf y$ is the globally asymptotically stable equilibrium of the subsystem on $U$. Moreover, $J\not\subseteq I$: otherwise $U=I$, the equilibrium $\mathbf x$ would be saturated for the same subsystem, and uniqueness would give $\mathbf x=\mathbf y$.

By \eqref{eq:mixed}, $r_j(\mathbf x)>0$ for every $j\in U\setminus I=J\setminus I$. This is the hypothesis of \cite[Theorem~18]{AKLS} for the equilibrium $\mathbf x$ and the set $U\supseteq I$. That theorem assumes VL stability only and hyperbolicity is not needed, see \cite[Remark~25]{AKLS}. It provides a solution $\boldsymbol\xi$ of \eqref{eq:lv} such that $\boldsymbol\xi(t)\to\mathbf x$ as $t\to-\infty$ and $\boldsymbol\xi(t)$ converges, as $t\to\infty$, to the globally asymptotically stable equilibrium of the subsystem on $U$, which is $\mathbf y$. This solution connects $\mathbf{x}$ with $\mathbf{y}$, so $\mathbf{x}\cto\mathbf{y}$.
\end{proof}

The above result implies that in the nonhyperbolic case the entry and exclusion conditions play assymetric role. New species cannot emerge from the source if the invasion rate is zero but it may disappear at a target with zero rate. 

\begin{remark}
For a VL stable matrix $A$ and any $\mathbf b$, Theorem~\ref{thm:exactcg} gives the following finite construction of the three graphs.
\begin{enumerate}[label={\bf (\arabic*)}]
	\item For every nonempty support $I$, solve $A_I\mathbf{x}_I=-\mathbf{b}_I$. Keep $\mathbf{x}_I$ if  $\mathbf{x}_I>\mathbf{0}$, and extend it by zero. Also include the empty support solution $\mathbf{0}$. Note tht all systems have unique solutions as VL stability implies that all principal submatrices are nonsingular.
	\item Compute every invasion rate $r_j(\mathbf{x})=b_j+(A\mathbf{x})_j$.
	\item In order to construct $\mathrm{CG}$ for each ordered pair of equilibria, draw the edge between them exactly
	when \eqref{eq:mixed} holds. If inequalities on both sides are strict, add the edge to $\mathrm{IG}$. If both are nonstrict add the edge to $\mathrm{WG}$.
\end{enumerate}

The construction requires the exact signs of the invasion rates: for rational $A$ and $\mathbf b$ they are provided by exact rational linear algebra, whereas floating-point values alone cannot certify a zero rate.
\end{remark}

\subsection{\texorpdfstring{Both inclusions in \eqref{eq:graphs} can be strict}{Both inclusions can be strict}}

To show $\mathrm{IG} \subsetneq \mathrm{CG}$, take
\[
 A=\begin{pmatrix}-1&-1\\0&-1\end{pmatrix},\qquad \mathbf{b}=\begin{pmatrix}1\\1\end{pmatrix},
 \qquad \mathbf{x}=\begin{pmatrix}1\\0\end{pmatrix},\quad \textrm{and}\quad  \mathbf{y}=\begin{pmatrix}0\\1\end{pmatrix},
\]
for which we have $r_2(\mathbf{x})=1$ and $r_1(\mathbf{y})=0$. Hence, by Theorem \ref{thm:exactcg}, $\mathbf{x} \cto \mathbf{y}$. However, the edge $\mathbf{x}\longrightarrow \mathbf{y}$ does not exist in the invasion graph $\mathrm{IG}$. An explicit connecting solution is 
\[
 \boldsymbol{\xi}(t)=\frac{e^t}{1+e^t} \begin{pmatrix}\frac{1}{\log(1+e^t)}\\ 1 \end{pmatrix} \quad \text{ for } t\in \mathbb{R},
\]
which is bounded and complete, with limits $\mathbf x$ and $\mathbf y$. Here $\rho(B)=0$, so $A\in\mathcal{H}_-\subseteq\mathcal{V}$.

To show that $\mathrm{CG}\subsetneq \mathrm{WG}$, consider
\[
 A=\begin{pmatrix}-1&0&0\\0&-1&0\\-1&1&-1\end{pmatrix},
 \qquad \mathbf{b} =\begin{pmatrix}1\\1\\1 \end{pmatrix}, \qquad \mathbf{x}=\begin{pmatrix}1\\0\\0 \end{pmatrix},\quad \textrm{and}\quad  \mathbf{y}=\begin{pmatrix}1\\1\\1 \end{pmatrix}.
\]
Again $\rho(B)=0$, so $A\in\mathcal{H}_-\subseteq\mathcal{V}$. Note that $\mathbf{x}\wto \mathbf{y}$, because $r_2(\mathbf{x})=1$ and $r_3(\mathbf{x})=0$, with no species lost. However, Proposition~\ref{lem:strictentry} shows that the edge $\mathbf{x}\cto\mathbf{y}$ does not exist, since this would imply that $r_3(\mathbf{x})>0$. Theorem~\ref{thm:exactcg} nevertheless gives the two-step path
\[
\mathbf{x}\cto \begin{pmatrix} 1\\1\\0 \end{pmatrix} \cto  \mathbf{y}.
\]
A path through intermediate equilibria need not give a direct connection, see also \cite[Example~34]{AKLS}.

\subsection{Sharpness of the acyclicity threshold}\label{sec:sharp}
Weak cycles can occur at $\rho(B)=1$. For example, take
\[
A=\begin{pmatrix}-1&-1\\-1&-1\end{pmatrix}, \qquad \mathbf{b}=\begin{pmatrix} 1\\ 1 \end{pmatrix}, \qquad \mathbf{x} = \begin{pmatrix} 1\\ 0 \end{pmatrix} \quad \text{ and } \quad \mathbf{y} = \begin{pmatrix} 0\\ 1 \end{pmatrix}.
\] 
We have $\rho(B)=1$ and $\mathbf{x}\wto \mathbf{y}\wto \mathbf{x}$, since $r(\mathbf{x})=r(\mathbf{y})=\mathbf 0$. It is a weak cycle but not a cycle in  $\mathrm{IG}$, since all rates vanish, nor in CG: $\frac{d}{dt}\log(u_1/u_2) = 0$ whenever both coordinates are positive, so no trajectory connects the two axes. Note that $A\not\in \mathcal{V}$ (as we will see later in Lemma~\ref{lem:energy} weak two-cycles in $\mathcal{V}$ are impossible for $N=2$). The matrix is singular and the whole segment $u_1+u_2=1$ consists of equilibria. 

\begin{proposition}\label{prop:strictboundary}
Assume that there is a cycle in $\mathrm{IG}$ for which the resident  interaction matrix of each equilibrium in this cycle is nonsingular. Then $\rho(B)>1$. In particular, if $A\in\mathcal{V}$ and there is a cycle in $\mathrm{IG}$, then $\rho(B)>1$.
\end{proposition}
\begin{proof}
If $\rho(B)\le1$, consider $A_\varepsilon=A+\varepsilon \diag(a_{11},a_{22},\ldots,a_{NN})$ and keep $\mathbf{b}$ fixed. The new matrix is $B_\varepsilon=\frac{1}{1+\varepsilon}B$, so $A_\varepsilon \in \mathcal{H}_-$. Since each $A_I$ is nonsingular, every equilibrium of the cycle can, for sufficiently small $\varepsilon>0$, be continues to the equilibrium  with the same support, and with the same strict inequalities that define the edges. This contradicts Corollary~\ref{cor:acyclic}, proving $\rho(B)>1$. If $A\in\mathcal{V}$, each $A_I$ is nonsingular, so the first assertion applies.
\end{proof}

For VL stable matrices, $\rho(B)\le1$ excludes cycles in $\mathrm{IG}$. The following May--Leonard type family shows what can happen in the connection graph: at $\rho(B)=1$ it has a heteroclinic cycle but no cycle in $\mathrm{IG}$, while above the threshold the same cycle can belong to both graphs. Here and below, $\mathbf{e}_i$ denotes the $i$-th vector of the standard basis of $\R^N$.

\begin{proposition}\label{th:sharp}
Let $c\ge0$, $\mathbf{b}=\one$, and $A_c=-\Id-C$ with 
\begin{equation}\label{eq:cyclic}
C=\begin{pmatrix}0&c&0\\0&0&c\\c&0&0\end{pmatrix}.
\end{equation}
The vectors $\mathbf{e}_1$, $\mathbf{e}_2$ and $\mathbf{e}_3$ are the single-species equilibria of \eqref{eq:lv}, and the following properties hold:
\begin{enumerate}[label={\bf (\roman*)}]
\item $\rho(B_c)=c$, so $A_c\in\mathcal{H}_-$ exactly for $0\le c<1$.
\item $A_c\in \mathcal{V}$ exactly for $0\le c<2$.
\item For $1<c<2$ we have
\[
\mathbf{e}_1 \longrightarrow \mathbf{e}_2 \longrightarrow \mathbf{e}_3 \longrightarrow \mathbf{e}_1
\quad \text{ and } \quad 
\mathbf{e}_1 \cto \mathbf{e}_2 \cto \mathbf{e}_3 \cto \mathbf{e}_1,
\]
and the same heteroclinic cycle exists at $c=1$, where its edges are weak edges but not edges of $\mathrm{IG}$.
\end{enumerate}
\end{proposition}
\begin{proof}
\textup{(i)} The matrix \eqref{eq:B} built from $A_c$ is $B_c=C$, and $C^3=c^3\Id$. Hence every eigenvalue of $B_c$ has modulus $c$.

\textup{(ii)} If $0\le c<2$, then taking $H=\Id$ we obtain
\[
A_c+A_c^{\mathsf T} = \begin{pmatrix} -2 & -c & -c \\ -c & -2 & -c \\ -c & -c & -2\end{pmatrix},
\]
which has eigenvalues $-2-2c$ and $-2+c$. Hence for $0\le c<2$ we have $A_c+A_c^{\mathsf T}\prec 0$ and $A_c\in \mathcal{V}$. For the converse, assume that $H=\diag(h_1,h_2,h_3)>0$ is such that $HA_c+A_c^{\mathsf T}H\prec 0$. We have
\[
HA_c+A_c^{\mathsf T}H = \begin{pmatrix} -2h_1 & -ch_1 & -ch_3 \\ -ch_1 & -2h_2 & -ch_2 \\ -ch_3 & -ch_2 & -2h_3\end{pmatrix}.
\]
The principal minors of order two give
\[
h_2 > \frac{c^2}{4}h_1, \qquad h_3 > \frac{c^2}{4}h_2, \quad \text{ and } \quad h_1 > \frac{c^2}{4}h_3.
\]
Multiplying these three inequalities we obtain $1>\frac{c^6}{64}$, that is, $c<2$.

\textup{(iii)} For $1<c<2$, the rates are $r_2(\mathbf{e}_1)=1$ and $r_1(\mathbf{e}_2)=1-c<0$; the other edges follow by cyclic permutation. At $c=1$, the invasion rates in the exclusion conditions vanish, so these are weak edges but not edges of $\mathrm{IG}$. On the face $u_3=0$, equation \eqref{eq:lv} becomes
\[
\left\{ \begin{aligned} 
& \dot u_1=u_1(1-u_1-cu_2),\\
& \dot u_2=u_2(1-u_2).
\end{aligned}\right.
\]
Take $u_2(t)=\frac{e^t}{1+e^t}$ and
\begin{equation}\label{eq:connection}
u_1(t)=\begin{cases}
\displaystyle\frac{(1-c)e^t}{(1+e^t)-(1+e^t)^c},&c\ne1,\\[7pt]
\displaystyle\frac{e^t}{(1+e^t)\log(1+e^t)},&c=1.
\end{cases}
\end{equation}
For $c\ge1$, the solution $\boldsymbol \xi(t)=(u_1(t),u_2(t),0)^{\mathsf T}$ is bounded, lies in the relative interior of the face $u_3=0$, and satisfies
\[
\mathbf{e}_1 \stackrel{t\to -\infty}{\longleftarrow}\boldsymbol{\xi}(t)\stackrel{t\to \infty}{\longrightarrow} \mathbf{e}_2.
\]
Cyclic permutation gives the other two connections.
\end{proof}

\subsection{The condition is not necessary for acyclicity}\label{sec:notnecessary}

A classical quadratic potential 
that goes back to MacArthur~\cite{MacArthur1969,MacArthur1970} gives acyclicity for symmetric Lotka--Volterra systems; see \cite[Section 7.2 and Lemma 37]{AKLS}. We state the argument for diagonally symmetrizable matrices and the weak graph. It applies for every $\mathbf b\in\mathbb R^N$, without a restriction on $\rho(B)$.

\begin{proposition}\label{prop:symmetric}
Assume there exists a positive diagonal matrix $H$ for which $HA$ is symmetric negative definite (thus $A\in \mathcal{V}$) and consider the map $\Phi\colon \mathbb{R}^N\to \mathbb{R}$ given by
\[
\Phi(\mathbf{u})=\mathbf{b}^{\mathsf T} H\mathbf{u} + \tfrac12 \mathbf{u}^{\mathsf T}HA\mathbf{u}. 
\]
Then every weak edge  $\mathbf{x}\wto \mathbf{y}$ satisfies $\Phi(\mathbf{x})<\Phi(\mathbf{y})$. Consequently, $\mathrm{WG}$ is acyclic. 
\end{proposition}
\begin{proof}
	Let $\mathbf{x}\wto\mathbf{y}$, put $\mathbf{z}=\mathbf{y}-\mathbf{x}$, and write $I=\supp\mathbf{x}$ and $J=\supp\mathbf{y}$. The resident equations give $r_j(\mathbf{y})=0$ for $j\in J$. For $i\in I\setminus J$, we have $z_i=-x_i<0$ and $r_i(\mathbf{y})\le0$. Outside $I\cup J$, $z_i=0$. Since $H$ is positive diagonal, $\mathbf{z}^{\mathsf T}Hr(\mathbf{y})\ge0$. 	We now use the exact Taylor formula for the quadratic function $\Phi$. Since $HA$ is symmetric, we have $\nabla \Phi(\mathbf y)=Hr(\mathbf y)$, and $\nabla^2 \Phi =HA$.
	Writing $\mathbf x=\mathbf y- \mathbf z$, we obtain
	\[
	\Phi(\mathbf x)
	=\Phi(\mathbf y)-\mathbf z^{\mathsf T}Hr(\mathbf y)
	+\tfrac12\mathbf z^{\mathsf T}HA\mathbf z.
	\]
Rearranging gives
\[
	\Phi(\mathbf{y})-\Phi(\mathbf{x})
	=\mathbf{z}^{\mathsf T}Hr(\mathbf{y})
	-\tfrac12\mathbf{z}^{\mathsf T}HA\mathbf{z}>0,
	\]
	because $\mathbf{z}\ne\mathbf 0$ and $HA\prec0$.
\end{proof}

Writing
\[
 \mathcal A_{\rm acyc}=\{A\in\mathcal{V} \colon \mathrm{WG} \text{ is acyclic for every }\mathbf{b}\in\R^N\},
\]
we therefore have, for $N\ge3$,
\begin{equation}\label{eq:acyclicclass}
 \mathcal{H}_-\subseteq\mathcal A_{\rm acyc}\subseteq \mathcal{V},
\end{equation}
with both inclusions strict. 

\section{Quantitative symmetric dominance}\label{sec:skew}

In Section~\ref{sec:dynamics} we showed that $A\in\mathcal{H}_-$ excludes cycles in the weak graph and hence in the invasion and connection graphs. By Proposition~\ref{prop:symmetric}, symmetry and negative definiteness of  $HA$, for some positive diagonal matrix $H$, also exclude weak cycles. In \cite{BKLM}, we asked whether dominance of the symmetric part also ensures the exclusion of cycles. The answer depends on how dominance is measured. If we measure dominance as the ratio of the entrywise norms of the skew and symmetric parts considered in \cite{BKLM}, the answer is negative, we show this in the example with matrices $A_\gamma$ below. We get the positive answer in Theorem~\ref{thm:skew} and condition~\eqref{eq:uniform} by using  the operator norm of the skew part $K$ weighted by the symmetric part $G$. The resulting criterion applies to matrices in $\mathcal{V}$ beyond both $\mathcal{H}_-$ and the diagonally symmetrizable case. 

For $A\in\mathcal{V}$ and a positive diagonal matrix $H=\diag(h_1,\ldots,h_N)$ such that $HA+A^{\mathsf T}H\prec0$, we define
$
G=-(HA+A^{\mathsf T}H)$, $K=HA-A^{\mathsf T}H$, $\Phi(\mathbf{x})=\tfrac12\mathbf{b}^{\mathsf T}H\mathbf{x}$, $\|\mathbf{z}\|_G^2=\mathbf{z}^{\mathsf T}G\mathbf{z}$, where $\mathbf{x},\mathbf{z}\in\R^N$. Thus $G$ is symmetric positive definite, $K$ is skew-symmetric, and $HA=\tfrac12(K-G)$.
At a nonnegative equilibrium $\mathbf{x}$ of \eqref{eq:lv}, the  equilibrium equations give $
 \Phi(\mathbf{x})=\tfrac14 \mathbf{x}^{\mathsf T}G\mathbf{x}$. Recall that $r(\mathbf{u})=\mathbf{b}+A\mathbf{u}$. At an equilibrium we have $\mathbf{x}^{\mathsf T}Hr(\mathbf{x})=0$, hence $\mathbf{x}^{\mathsf T}HA\mathbf{x}=-\mathbf{b}^{\mathsf T}H\mathbf{x}$, so that $\Phi(\mathbf{x})$ is also the value at $\mathbf{x}$ of the function $\Phi$ of Proposition~\ref{prop:symmetric}.
 
When $HA$ is not symmetric, the change of $\Phi$ along a weak edge has a symmetric and a skew contribution, which we now separate.

\begin{lemma}\label{lem:energy}
Let $A\in\mathcal{V}$, and let $H$, $G$, $K$ and $\Phi$ be as above. For a weak edge $\mathbf{x}\wto\mathbf{y}$ with supports $I=\supp \mathbf{x}$ and $J=\supp \mathbf{y}$, define the two weighted total invasion rates
\[
 T_1=\sum_{j\in J\setminus I}h_jr_j(\mathbf{x})y_j\ge0,\qquad \textrm{and}\qquad 
 T_2=\sum_{i\in I\setminus J}h_ir_i(\mathbf{y})x_i\le0.
\]
Then
$
 T_1+T_2=\tfrac12\norm{\mathbf{y}-\mathbf{x}}_G^2,$ $T_1-T_2\ge\tfrac12\norm{\mathbf{y}-\mathbf{x}}_G^2$, and $
 \Phi(\mathbf{y})-\Phi(\mathbf{x})=\tfrac12(T_1-T_2)+\tfrac12 \mathbf{x}^{\mathsf T}K\mathbf{y}$.
Consequently, every weak edge introduces at least one species with a positive invasion rate, so an edge with $\supp\mathbf{y}\subseteq\supp\mathbf{x}$ is impossible. Moreover, there are no weak two-cycles.
\end{lemma}
\begin{proof}
Since $r_i(\mathbf{x})=0$ for every $i\in\supp\mathbf{x}$, we obtain
\[
\begin{aligned}
r(\mathbf{x})^{\mathsf T}H(\mathbf{y}-\mathbf{x})
&=\sum_{j=1}^N r_j(\mathbf{x})h_j(y_j-x_j)=\sum_{j\not\in I}r_j(\mathbf{x})h_j(y_j-x_j)=T_1,
\end{aligned}
\]
and, analogously, 
$ r(\mathbf{y})^{\mathsf T}H(\mathbf{y}-\mathbf{x})=-T_2.
$
Subtraction gives
$-T_1-T_2=(\mathbf{y}-\mathbf{x})^{\mathsf T}HA(\mathbf{y}-\mathbf{x})=-\norm{\mathbf{y}-\mathbf{x}}_G^2/2$.
Also $\mathbf{y}^{\mathsf T}Hr(\mathbf{x})=T_1$ and $\mathbf{x}^{\mathsf T}Hr(\mathbf{y})=T_2$, hence
\[
\begin{aligned}
T_1-T_2
&=\mathbf b^{\mathsf T}H(\mathbf{y}-\mathbf{x})
  +\mathbf y^{\mathsf T}HA\mathbf{x}-\mathbf{x}^{\mathsf T}HA\mathbf{y}=2\bigl(\Phi(\mathbf{y})-\Phi(\mathbf{x})\bigr)-\mathbf{x}^{\mathsf T}K\mathbf{y}.
\end{aligned}
\]
The desired inequality follows from 
\[
T_1-T_2 = (T_1+T_2)-2T_2 = \frac12\|\mathbf{y}-\mathbf{x}\|_G^2 - 2T_2 \ge \frac12\|\mathbf{y}-\mathbf{x}\|_G^2,
\]
since $T_2\le0$.
The identity $T_1=\frac12\|\mathbf{x}-\mathbf{y}\|_G^2-T_2$ implies that $T_1>0$ for distinct endpoints $\mathbf{x}$ and  $\mathbf{y}$. Hence some $j\in J\setminus I$ satisfies $r_j(\mathbf{x})>0$.  

To exclude a cycle with two equilibria, suppose that we have the weak edges  $\mathbf{x}\wto\mathbf{y}\wto\mathbf{x}$. For the edge $\mathbf{y}\wto\mathbf{x}$ the two sums that define $T_1$ and $T_2$ exchange their roles. We denote them by $\tilde T_1$ and $\tilde T_2$ and we have $\tilde{T}_1=T_2$ and $\tilde{T}_2=T_1$. The sign conditions $\tilde{T}_1\ge0\ge T_2$ and $T_1\ge0\ge\tilde{T}_2$ therefore give $T_1=T_2=0$. Then $\tfrac12\norm{\mathbf{y}-\mathbf{x}}_G^2=T_1+T_2=0$, which contradicts $\mathbf{x}\ne\mathbf{y}$.
\end{proof}

If $HA$ is symmetric, then $K=0$, and Lemma~\ref{lem:energy} gives $\Phi(\mathbf{y})-\Phi(\mathbf{x})=\tfrac12(T_1-T_2)\ge\tfrac14\norm{\mathbf{y}-\mathbf{x}}_G^2>0$ along every weak edge $\mathbf{x}\wto\mathbf{y}$. This recovers Proposition~\ref{prop:symmetric}, with a quantitative lower bound.

To control the skew contribution around a cycle we need an estimate for cyclic sequences of vectors. It rests on the following discrete Wirtinger inequality of Fan, Taussky and Todd~\cite{FanTausskyTodd1955}.
\begin{lemma}[Fan--Taussky--Todd {\cite[Theorem~10, p.~80]{FanTausskyTodd1955}}]\label{lemma:isoper}
	Let $m\ge3$ and let $a_1,\ldots,a_m\in\mathbb R$ satisfy
	$\sum_{k=1}^m a_k=0$, with $a_{m+1}=a_1$. Then
	\begin{equation}\label{eq:wirtinger}
		\sum_{k=1}^m(a_{k+1}-a_k)^2
		\ge 4\sin^2\!\left(\frac{\pi}{m}\right)\sum_{k=1}^m a_k^2.
	\end{equation}
\end{lemma}
The next lemma is a matrix form of the discrete isoperimetric inequality. Its planar specialization bounds the absolute signed area of a polygon by  a multiple of the sum of squared side lengths. 
\begin{lemma}\label{lem:skew-area}
	Let $L\in\mathbb R^{N\times N}$ be a skew-symmetric matrix and let
	$\mathbf z_1,\ldots,\mathbf z_m\in\mathbb R^N$ be arbitrary vectors, where
	$m\ge3$ and $\mathbf z_{m+1}=\mathbf z_1$. Then
	\begin{equation}\label{eq:skew}
		\left|\sum_{k=1}^m\mathbf z_k^{\mathsf T}L\mathbf z_{k+1}\right|
		\le\frac{\|L\|_2}{2}\cot\!\left(\frac{\pi}{m}\right)
		\sum_{k=1}^m\|\mathbf z_{k+1}-\mathbf z_k\|_2^2,
	\end{equation}
	where $\|L\|_2$ is the Euclidean operator norm.
\end{lemma}

\begin{proof}
	Both sides of \eqref{eq:skew} do not change if we substract a given vector from every $\mathbf z_i$.
	The skew symmetric matrix on the left makes all  additional
	terms in the cyclic sum cancel after such translation. We subtract from all vectors their arithmetic  mean.
	This means that without loss of generality we may assume $\sum_k\mathbf z_k=\mathbf 0$.
	Set $\theta=\pi/m$ and $E=\sum_k\|\mathbf z_{k+1}-\mathbf z_k\|_2^2$.
	Applying \eqref{eq:wirtinger} to each coordinate and summing gives
	\[
	4\sin^2\theta\sum_k\|\mathbf z_k\|_2^2\le E.
	\]
	Together with the parallelogram identity, this yields
	\[
	\sum_k\|\mathbf z_{k+1}+\mathbf z_k\|_2^2
	=4\sum_k\|\mathbf z_k\|_2^2-E
	\le\cot^2\theta\,E.
	\]
	Finally, skew symmetry and Cauchy--Schwarz give
	\begin{align*}
		\left|\sum_k\mathbf z_k^{\mathsf T}L\mathbf z_{k+1}\right|
		&=\frac12\left|\sum_k
		(\mathbf z_k+\mathbf z_{k+1})^{\mathsf T}
		L(\mathbf z_{k+1}-\mathbf z_k)\right|\\
		&\le\frac{\|L\|_2}{2}
		\left(\sum_k\|\mathbf z_{k+1}+\mathbf z_k\|_2^2\right)^{1/2}
		E^{1/2}\le\frac{\|L\|_2}{2}\cot\theta\,E.\qedhere
	\end{align*}
\end{proof}
We can now prove the quantitative criterion for excluding cycles.
\begin{theorem}\label{thm:skew}
	Let $A\in\mathcal V$, let $H$ be a positive diagonal matrix such that $G=-(HA+A^{\mathsf T}H)\succ0$, and put $K=HA-A^{\mathsf T}H$. If, for some $\mathbf b\in\R^N$, the weak graph of \eqref{eq:lv} contains a cycle through $m\ge3$ distinct equilibria, then
	\begin{equation}\label{eq:kappa}
		\norm{G^{-1/2}KG^{-1/2}}_2\ge\tan\!\left(\frac{\pi}{m}\right).
	\end{equation}
	Consequently, if $\norm{G^{-1/2}KG^{-1/2}}_2<\tan(\pi/m)$, then for every $\mathbf b\in\R^N$ the weak graph has no cycle of length at most $m$.
\end{theorem}
\begin{proof}
	First note that 	
Lemma~\ref{lem:energy} excludes cycles of length two. Let $\mathbf{x}_1,\ldots,\mathbf{x}_m$ be the equilibria of the cycle, in order, and put $\mathbf{x}_{m+1}=\mathbf{x}_1$. All sums below run over $k=1,\ldots,m$. Set $\mathbf z_k=G^{1/2}\mathbf x_k$ and
	\[
	E=\sum_k\norm{\mathbf z_{k+1}-\mathbf z_k}_2^2 .
	\]
	Since the equilibria are distinct and $G\succ0$, we have $E>0$. For the edge $\mathbf x_k\wto\mathbf x_{k+1}$, let $T_1^{(k)}$ and $T_2^{(k)}$ be the two sums of Lemma~\ref{lem:energy}. That lemma gives
	\[
	\Phi(\mathbf x_{k+1})-\Phi(\mathbf x_k)=\tfrac12\bigl(T_1^{(k)}-T_2^{(k)}\bigr)+\tfrac12\mathbf x_k^{\mathsf T}K\mathbf x_{k+1}
	\qquad\text{and}\qquad
	T_1^{(k)}-T_2^{(k)}\ge\tfrac12\norm{\mathbf x_{k+1}-\mathbf x_k}_G^2 .
	\]
	Around the cycle the left-hand sides of the identities add up to zero. Hence
	\[
	\sum_k \mathbf{x}_k^{\mathsf T}K\mathbf{x}_{k+1}
	=-\sum_k\bigl(T_1^{(k)}-T_2^{(k)}\bigr)\le -\frac{E}{2}.
	\]
	On the other hand, $L=G^{-1/2}KG^{-1/2}$ is skew-symmetric and $\mathbf x_k^{\mathsf T}K\mathbf x_{k+1}=\mathbf z_k^{\mathsf T}L\mathbf z_{k+1}$, so Lemma~\ref{lem:skew-area} gives
	\[
	\left|\sum_k\mathbf{x}_k^{\mathsf T}K\mathbf{x}_{k+1}\right|
	=\left|\sum_k\mathbf{z}_k^{\mathsf T}L\mathbf{z}_{k+1}\right|
	\le\frac{\norm{L}_2}{2}\cot\!\left(\frac{\pi}{m}\right)E .
	\]
	The two estimates together give $E/2\le\tfrac12\norm{L}_2\cot(\pi/m)\,E$. Dividing by $E/2>0$ yields \eqref{eq:kappa}.
\end{proof}

Shorter cycles thus require a larger skew part $K$ relative to $G$. To exclude cycles of every length, it remains to bound the length of a weak cycle.

\begin{lemma}\label{lem:cyclevertices}
	Let $A\in\mathcal V$ and $\mathbf b\in\R^N$. A weak cycle contains neither the zero equilibrium, nor an equilibrium with $N$ species, nor an equilibrium with $N-1$ species. Consequently, it has at most $2^N-N-2$ vertices.
\end{lemma}
\begin{proof}
	By Lemma~\ref{lem:energy}, for every weak edge $\mathbf x\wto\mathbf y$ there is some $j\in\supp\mathbf y\setminus\supp\mathbf x$ with $r_j(\mathbf x)>0$. Hence no weak edge ends at $\mathbf 0$ and no weak edge starts at an equilibrium with $N$ species, so neither of them lies on a weak cycle.
	
	Now let $\mathbf x\wto\mathbf y$ with $I=\supp\mathbf x$, $|I|=N-1$, and $J=\supp\mathbf y$. Since $J\setminus I\ne\varnothing$, we have $I\cup J=\{1,\ldots,N\}$. Every $i\notin J$ therefore belongs to $I\setminus J$, and \eqref{eq:weak} gives $r_i(\mathbf y)\le0$. So no absent species has a positive invasion rate at $\mathbf y$, and by the first paragraph no weak edge starts at $\mathbf y$. Thus $\mathbf x$ lies on no weak cycle.
	
	Finally, there are $2^N$ supports, each carrying at most one equilibrium because $A\in\mathcal V$, and $N+2$ of them have just been excluded.
\end{proof}
\begin{corollary}\label{cor:uniform}
	Let $N\ge3$ and $A\in\mathcal V$, and let $H$, $G$ and $K$ be as in Theorem~\ref{thm:skew}. If
	\begin{equation}\label{eq:uniform}
		\norm{G^{-1/2}KG^{-1/2}}_2 <\tan\!\left(\frac{\pi}{2^N-N-2}\right),
	\end{equation}
	then the weak graph is acyclic for every $\mathbf b\in\R^N$. Consequently, the invasion and connection graphs are acyclic, and the semiflow $\mathcal T$ is gradient with respect to its equilibria.
\end{corollary}

\begin{proof}
	By Lemma~\ref{lem:cyclevertices}, a weak cycle has at most $2^N-N-2$ vertices, and this number is at least $3$. Theorem~\ref{thm:skew} with $m=2^N-N-2$ excludes all weak cycles of at most this length. The remaining assertions follow from \eqref{eq:graphs} and Corollary~\ref{cor:wggradient}.
\end{proof}

The threshold in \eqref{eq:uniform} decays exponentially with the number of species: $\tan\bigl(\pi/(2^N-N-2)\bigr)\sim\pi\,2^{-N}$ as $N\to\infty$. For large $N$ the uniform criterion is therefore very restrictive, and Theorem~\ref{thm:skew}, which excludes cycles up to a given length, is the more useful statement.

Note that Theorem~\ref{thm:skew} can alternatively be derived from the characterization of cyclically monotone linear operators in \cite[Theorem~1]{Shiu1976}: the condition $\norm{G^{-1/2}KG^{-1/2}}_2\le\tan(\pi/m)$ means that the linear map $-HA=\tfrac12(G-K)$ is $m$-cyclically monotone. 

The quantity $\norm{G^{-1/2}KG^{-1/2}}_2$ has the variational characterization 
\[
\norm{G^{-1/2}KG^{-1/2}}_2 = \max_{\mathbf{x},\mathbf{y}\ne\mathbf 0}\frac{|\mathbf{x}^{\mathsf T}K\mathbf{y}|}{\sqrt{\mathbf{x}^{\mathsf T}G\mathbf{x}}\sqrt{\mathbf{y}^{\mathsf T}G\mathbf{y}}},
\]
from which it follows that 
\[
\norm{G^{-1/2}KG^{-1/2}}_2 \le \frac{\|K\|_2}{\lambda_{\min}(G)}. 
\]
Therefore the condition
$
\frac{\|K\|_2}{\lambda_{\min}(G)} < \tan(\pi/m)
$ excludes cycles of length at most $m$, and \[
\frac{\|K\|_2}{\lambda_{\min}(G)} < \tan\left(\frac{\pi}{2^N-N-2}\right)
\] excludes all cycles.
When $A=A_{\mathrm{sym}}+A_{\mathrm{skew}}$, $A_{\mathrm{sym}}=A_{\mathrm{sym}}^{\mathsf T}\prec0$ and $A_{\mathrm{skew}}=-A_{\mathrm{skew}}^{\mathsf T}$, one may use $H=\Id$. Then $G=-2A_{\mathrm{sym}}$ and $K=2A_{\mathrm{skew}}$, so  $\norm{G^{-1/2}KG^{-1/2}}_2=\norm{(-A_{\mathrm{sym}})^{-1/2}A_{\mathrm{skew}}(-A_{\mathrm{sym}})^{-1/2}}_2$. Condition
\eqref{eq:uniform} then reads
\[
 \norm{(-A_{\mathrm{sym}})^{-1/2}A_{\mathrm{skew}}(-A_{\mathrm{sym}})^{-1/2}}_2
       <\tan\!\left(\frac{\pi}{2^N-N-2}\right).
\]
This gives an explicit sufficient condition for  cycle exclusion in terms of the symmetric and skew parts, answering Open Question~2 of \cite{BKLM}. There we considered the ratio of the entrywise norms of the skew and symmetric parts, where $\norm{C}_{\mathrm{ent}}=\sum_{i,j=1}^N|c_{ij}|$. This ratio alone is insufficient to exclude cycles. Indeed, for $\gamma>0$, consider 
\[
 A_\gamma =\begin{pmatrix}
 -1&-2&2&0\\2&-1&-2&0\\-2&2&-1&0\\0&0&0&-\gamma
 \end{pmatrix},\qquad \textrm{and}\qquad \mathbf{b}=\one.
\]
Here $H=\Id$, $A_\gamma+A_\gamma^{\mathsf T}=\diag(-2,-2,-2,-2\gamma)\prec0$, and the singleton equilibria $\mathbf{e}_1,\mathbf{e}_2,\mathbf{e}_3$ form the cycle $\{1\}\to \{2\} \to \{3\}\to\{1\}$ in $\mathrm{IG}$. Since $A_\gamma\in\mathcal V$, these edges are heteroclinic connections by Theorem~\ref{thm:exactcg}. Nevertheless,
\[
 \frac{\norm{(A_\gamma-A_\gamma^{\mathsf T})/2}_{\mathrm{ent}}}
      {\norm{(A_\gamma+A_\gamma^{\mathsf T})/2}_{\mathrm{ent}}}
       =\frac{12}{3+\gamma}\longrightarrow0 \quad \text{ as } \gamma\to \infty.
\]
Thus no positive threshold for this ratio alone can exclude cycles, even within $\mathcal{V}$.  Increasing the self-damping of the independent fourth species makes the ratio arbitrarily small while leaving the cycling subsystem unchanged. For $H=\Id$, in contrast, $\norm{G^{-1/2}KG^{-1/2}}_2=2\sqrt3$ for every $\gamma$.

The following example  satisfies the new criterion but neither generalized strict diagonal dominance nor diagonal symmetrizability. Consider the matrix
\[
A=\begin{pmatrix}
	-20&-9&-10&-11\\
	-11&-20&-9&-10\\
	-10&-11&-20&-9\\
	-9&-10&-11&-20
\end{pmatrix}.
\]
 For this matrix $\rho(B)=3/2>1$ and $A\notin\mathcal{H}_-$. Moreover, $HA$ cannot be made symmetric for a positive diagonal matrix $H$. But
\[
\|G^{-1/2}KG^{-1/2}\|_2=\frac15
<\tan\!\left(\frac{\pi}{2^4-4-2}\right)\approx0.324920,
\]
so the new criterion excludes all weak cycles for every $\mathbf b\in\R^4$. 

\appendix

\section{Characterization and properties of the class \texorpdfstring{$\mathcal{H}_-$}{H-}}\label{sec:equiv}

We collect the classical characterizations of generalized strict diagonal dominance and explain their relation to $\mathcal H_-$, $\mathcal V$, and ordinary row or column dominance. See \cite{FiedlerPtak1962,Plemmons1977,BermanPlemmons1994} and \cite[Section~4 and Table~1]{BCGM}. We include the proofs for completeness and claim no novelty in this Appendix.

We use the notation of Sections~\ref{sec:1} and~\ref{sec:dynamics}: the diagonal of $A$ is negative, $B$ is the matrix \eqref{eq:B}, and $M=\diag(|a_{ii}|)(\Id-B)$ is the comparison matrix of $A$.
All vector and matrix inequalities are understood entrywise.
A real square matrix is \emph{stable} if all its
eigenvalues have strictly negative real parts. A \emph{nonsingular $M$-matrix} is
a real square matrix with nonpositive off-diagonal entries whose inverse
exists and is entrywise nonnegative.

\begin{proposition}\label{th:equiv}
The following conditions are equivalent:
\begin{enumerate}[label={\bf (\roman*)},itemsep=1pt]
\item $\rho(B)<1$.
\item There exists $\mathbf v>\mathbf 0$ such that $M\mathbf v>\mathbf 0$, or equivalently
\begin{equation}\label{eq:weight}
 |a_{ii}|v_i>\sum_{j\ne i}|a_{ij}|v_j
 \qquad\text{for}\qquad i=1,\ldots,N.
\end{equation}
\item There exists $\mathbf w>\mathbf 0$ such that $M^{\mathsf T}\mathbf w>\mathbf 0$, or equivalently 
\begin{equation}\label{eq:col_weight}
	|a_{jj}|w_j > \sum_{i\ne j}|a_{ij}|w_i  \qquad\text{for}\qquad j=1,\ldots,N.
\end{equation}	
\item $M$ is invertible and $M^{-1}\ge0$.
\item Every principal minor of $M$ is positive.
\item The matrix $-M$ is stable.
\end{enumerate}
Conditions \textup{(ii)} and \textup{(iii)} express generalized strict
diagonal dominance by rows and by columns, respectively. In particular,
every matrix with negative diagonal that is strictly diagonally dominant
by rows or by columns belongs to $\mathcal H_-$. 

If these conditions hold, then $\mathbf v=M^{-1}\one$ and $\mathbf w=M^{-\mathsf T}\one$ are strictly positive,
and they are admissible weights in \textup{(ii)} and \textup{(iii)}, respectively.
\end{proposition}
\begin{proof}
If \textup{(i)} holds, the Neumann series gives
\[
 M^{-1}=(\Id-B)^{-1}\diag(|a_{ii}|)^{-1}
       =\sum_{k=0}^{\infty}B^k\diag(|a_{ii}|)^{-1}\ge \diag(|a_{ii}|)^{-1}.
\]
This proves \textup{(iv)}. Now assume \textup{(iv)}. Every row
and column of $M^{-1}$ has a positive entry, since $M^{-1}$ is invertible.
Consequently,
$
 \mathbf v=M^{-1}\one>\mathbf 0$ and $\mathbf w=M^{-\mathsf T}\one>\mathbf 0$. Since $M\mathbf v=M^{\mathsf T}\mathbf w=\one$,
conditions \textup{(ii)} and \textup{(iii)} hold. If \textup{(ii)}
holds, then $B\mathbf v<\mathbf v$. 
Since $B\ge0$, the maximum absolute row-sum norm gives
\[
 \rho(B)=\rho(\diag(\mathbf{v})^{-1}B\diag(\mathbf{v}))
 \le\|\diag(\mathbf{v})^{-1}B\diag(\mathbf{v})\|_\infty
 =\max_i\frac{(B\mathbf v)_i}{v_i}<1.
\]
This proves \textup{(i)}.  Condition \textup{(iii)} for $A$ is condition
\textup{(ii)} for $A^{\mathsf T}$: the comparison matrix of $A^{\mathsf T}$ is
$M^{\mathsf T}$, and the matrix \eqref{eq:B} built from $A^{\mathsf T}$ is
$\diag(|a_{ii}|)^{-1}B^{\mathsf T}\diag(|a_{ii}|)$, which has the same spectral radius as $B$. Hence
\textup{(iii)} implies \textup{(i)} by the argument just given. The equivalences 
$\textup{(iv)}\Leftrightarrow\textup{(v)}\Leftrightarrow\textup{(vi)}$
are classical characterizations of nonsingular $M$-matrices, see
\cite[Theorem~6.2.3]{BermanPlemmons1994} and
\cite{FiedlerPtak1962,Plemmons1977}.
\end{proof}

We next prove that matrices in $\mathcal H_-$ are Volterra--Lyapunov
stable by constructing a positive diagonal matrix $H$.

\begin{proposition}\label{lem:vl}
If $A\in\mathcal H_-$, then $A\in\mathcal V$. More precisely,
in the definition of VL stability one may take
\[
H=\operatorname{diag}\!\left(\frac{w_i}{v_i}\right)
\quad \text{where} \quad
\mathbf v=M^{-1}\mathbf1\quad \text{and}\quad 
\mathbf w=M^{-\mathsf T}\mathbf1.
\]
\end{proposition}
\begin{proof}
By Proposition~\ref{th:equiv}, the vectors
	$\mathbf v=M^{-1}\mathbf1$ and $\mathbf w=M^{-\mathsf T}\mathbf1$
	are strictly positive and satisfy $M\mathbf v=M^{\mathsf T}\mathbf w=\mathbf1$.
	Hence $H$ is positive diagonal. These identities give
	\[
	\sum_{j\ne i}|a_{ij}|v_j=-a_{ii}v_i-1
	\qquad \text{and}\qquad 
	\sum_{j\ne i}|a_{ji}|w_j=-a_{ii}w_i-1.
	\]
	We define  $
	R=\diag(\mathbf v)(HA+A^{\mathsf T}H)\diag(\mathbf v)$,
	a symmetric matrix with entries
	$R_{ii}=2a_{ii}v_iw_i$ and $
	R_{ij}=w_iv_ja_{ij}+v_iw_ja_{ji}$ for $i\ne j$.
	Consequently,
\begin{align*}
	\sum_{j\ne i}|R_{ij}|
	&\le w_i\sum_{j\ne i}|a_{ij}|v_j
	     +v_i\sum_{j\ne i}|a_{ji}|w_j\\
	&=w_i(-a_{ii}v_i-1)+v_i(-a_{ii}w_i-1)=-2a_{ii}v_iw_i-w_i-v_i<|R_{ii}|.
	\end{align*}
	Thus $R$ is symmetric, has negative diagonal, and is strictly diagonally
	dominant. It is therefore negative definite, and so is $HA+A^{\mathsf T}H=\diag(\mathbf v)^{-1}R\diag(\mathbf v)^{-1}$.
\end{proof}

Propositions~\ref{th:equiv} and~\ref{lem:vl} give the implications from ordinary row or column dominance to
generalized dominance, and from generalized dominance to VL stability.
Neither converse holds when $N\ge2$. Indeed,
generalized dominance need not imply ordinary row or column dominance.
For example, take
\[
 A=\begin{pmatrix}-1&2\\1/5&-1\end{pmatrix}.
\]
Its first row and second column fail strict diagonal dominance, while
$\rho(B)=\sqrt{2/5}<1$. The gap between ordinary and weighted dominance can be arbitrarily large. To see
this, fix $0<q<1$ and $\gamma>0$, and let
\[
A=\begin{pmatrix}-1&\gamma&0\\0&-1&\gamma\\q^3/\gamma^2&0&-1\end{pmatrix},\qquad B=\begin{pmatrix}0&\gamma&0\\0&0&\gamma\\q^3/\gamma^2&0&0\end{pmatrix}.
\]
We have $B^3=q^3\Id$, and hence $\rho(B)=q$. A possible choice of weights in condition \textup{(ii)} of Proposition~\ref{th:equiv} is $
\mathbf v=(\gamma^2/q^2,\gamma/q,1)^{\mathsf T}
$. For $\gamma>1$, ordinary row and
column dominance both fail, and their defects grow without bound as
$\gamma\to\infty$. 

 VL stability need not imply generalized
dominance either. Take
\[
 A=\begin{pmatrix}-1&\beta\\-\beta&-1\end{pmatrix},\qquad \beta>1.
\]
Here $A+A^{\mathsf T}=-2\Id$, so $H=\Id$ proves VL stability, whereas
$\rho(B)=\beta>1$. The spectral radius $\rho(B)$ depends only on the absolute values $|a_{ij}|$ and in this example $a_{12}$ and $a_{21}$ have opposite signs, and this is why $A+A^{\mathsf T}$ is  diagonal.


\end{document}